\documentclass[12pt]{amsart}
\usepackage{amscd}
\usepackage{pstcol,pst-plot,pst-3d}
\def\NZQ{\Bbb}               % the font for N,Z,Q,R,C
\def\NN{{\NZQ N}}
\def\QQ{{\NZQ Q}}
\def\ZZ{{\NZQ Z}}
\def\RR{{\NZQ R}}

\def\D{{\Delta}}
\def\G{{\Gamma}}
\def\mm{{\frak m}}
\def\F{{\mathcal F}}
\def\a{{\bold a}}
\def\b{{\bold b}}
\def\c{{\bold c}}

\def\x{{\bold x}}
\def\y{{\bold y}}
\def\z{{\bold z}}
\def\e{{\bold e}}
\def\1{{\mathbf 1}}
\def\0{{\mathbf 0}}
\def\opn#1#2{\def#1{\operatorname{#2}}} % to make operators
\opn\diam{diam} 
\opn\depth{depth}
\opn\star{st}
\opn\lk{lk}

\newtheorem{Theorem}{Theorem}[section]
\newtheorem{Lemma}[Theorem]{Lemma}
\newtheorem{Corollary}[Theorem]{Corollary}

\newtheorem{Remark}[Theorem]{Remark}

\newtheorem{Example}[Theorem]{Example}

\newtheorem{Question}[Theorem]{Question}
\begin{document}
\title{Cohen-Macaulayness of monomial ideals\\
and symbolic powers of Stanley-Reisner ideals}

\author{Nguyen Cong Minh}
\address{Department of Mathematics, University of Education, 136 Xuan Thuy, Hanoi, Vietnam}
\email{ngcminh@gmail.com}

\author{Ngo Viet Trung}
\address{Institute of Mathematics, 18 Hoang Quoc Viet, Hanoi, Vietnam}
\email{nvtrung@math.ac.vn}
\keywords{Cohen-Macaulayness, monomial ideal, linear inequalities, simplicial complex, Stanley-Reisner ideal, symbolic power, graph, matroid complex}
\thanks{The authors are supported by the National Foundation of  
Science and Technology Development} 

\maketitle

\begin{abstract}
We present criteria for the Cohen-Macaulayness of a monomial ideal in terms of its primary decomposition.
These criteria allow us to use tools of  graph theory and of linear programming to study the Cohen-Macaulayness of
monomial ideals which are intersections of prime ideal powers. We can characterize the Cohen-Macaulayness of the second symbolic power or of all symbolic powers of a Stanley-Reisner ideal in terms of the simplicial complex. These characterizations show that the simplicial complex must be very compact if some symbolic power is Cohen-Macaulay. In particular,  all symbolic powers are Cohen-Macaulay  if and only if the simplicial complex is a matroid complex. We also prove that the Cohen-Macaulayness can pass from a symbolic power to another symbolic powers in different ways. 
\end{abstract}

\section*{Introduction}

The main aim of this paper is to characterize the Cohen-Macaulayness of symbolic powers of a squarefree monomial ideal in terms of the associated simplicial complex. This problem arises when we want to study the Cohen-Macaulayness of ordinary powers of a squarefree monomial ideal. Recall that the $m$-th symbolic power $I^{(m)}$ of an ideal $I$ in a Noetherian ring is defined as the intersection of the primary components of $I^m$ associated with the minimal primes.
For a radical ideal in a polynomial ring over a field of characteristic zero, Nagata and Zariski showed that $I^{(m)}$ is the ideal of the polynomials that vanish to order $m$ on the affine variety $V(I)$.
The usual way for testing the Cohen-Macaulayness of a monomial ideal is to pass to the polarized ideal in order to apply Reisner's criterion for squarefree monomial ideals. To polarize an ideal  we have to know the generators, which are not available for symbolic powers. So we need to find necessary and sufficient conditions for a monomial ideal to be Cohen-Macaulay in terms of its primary decomposition.
\smallskip

Recently Takayama \cite{Ta} gave a formula for the local cohomology modules of an arbitrary monomial ideal by means of certain simplicial complexes associated with each degree of the multigrading. The formula is technically complicated and involves the generators of the ideal. In \cite{MT} we succeeded in using Takayama's formula to characterize the Cohen-Macaulayness of symbolic powers of two-dimensional squarefree monomial ideals. Inspired of \cite{MT} we shall show in Section 1 that Takayama's formula actually yields the following criteria for the Cohen-Macaulayness of a monomial ideal in terms of its primary decomposition.
\smallskip

Let $I$ be a monomial ideal in the polynomial ring $S = k[x_1,...,x_n]$, where $k$ is a field of arbitrary characteristic. 
Let $\D$ be the simplicial complex on $[n] =\{1,...,n\}$ such that $\sqrt{I}$ is the Stanley-Reisner ideal  
$$I_\D = \bigcap_{F \in \F(\D)}P_F,$$
where $\F(\D)$ denotes the set of the facets of $\D$ and $P_F$ is the prime ideal of $S$ generated by the variables $x_i$, $i \not\in F$. Assume that 
$$I = \bigcap_{G \in \F(\D)}I_F,$$
where $I_F$ is the $P_F$-primary component of $I$.
\medskip

For every point $\a  = (a_1,...,a_n) \in \NN^n$ we set $x^\a = x_1^{a_1}\cdots x_n^{a_n}$ and we denote by $\D_\a$ the simplicial complex on $[n]$ with $\F(\D_\a) = \{F \in \F(\D)|\ x^\a \not\in I_F\}.$ Moreover, for every simplicial complex $\G$ with $\F(\G) \subseteq \F(\D)$ we set
$$L_\G(I) := \big\{\a \in \NN^n \big|\  x^\a \in \bigcap_{F \in \F(\D)\setminus\F(\G)}I_F \setminus \bigcup_{G \in \F(\G)}I_G\big\}.$$
\smallskip

\noindent {\bf Theorem 1.6.} 
Assume that $I$ is an unmixed monomial ideal. Then the following conditions are equivalent:\par
{\rm (i)} $I$ is a Cohen-Macaulay ideal,\par
{\rm (ii)} $\D_\a$ is a Cohen-Macaulay complex for all $\a \in \NN^n$,\par
{\rm (iii)} $L_\G(I) = \emptyset$ for every non-Cohen-Macaulay complex $\G$ with $\F(\G) \subseteq \F(\D)$.
\medskip

Here we call a simplicial complex $\G$ Cohen-Macaulay if $\tilde H_j(\lk_\G F,k) = 0$ for all $F \in \G$, $j < \dim \lk_\G F$. We can easily deduce from Theorem \ref{criteria}(ii) previous results on the Cohen-Macaulayness of squarefree monomial ideals such as Reisner's criterion that $I_\D$ is Cohen-Macaulay if and only if $\D$ is Cohen-Macaulay \cite{Re} and Eisenbud's observation that $\sqrt{I}$ is Cohen-Macaulay if $I$ is Cohen-Macaulay \cite{HTT}.
\smallskip

Theorem \ref{criteria}(iii) is especially useful when $I$ is the intersection of prime ideal powers, that is, all primary components $I_F$ are of the form  $P_F^m$ for some positive integers $m$. In this case, $x^\a \in I_F$ if and only if $\sum_{i\not\in F}a_i \le m$. Hence, $L_\G(I)$ is the set of solutions in $\NN^n$ of a system of linear inequalities. So we only need to test the inconsistency of systems of linear inequalities associated with the non-Cohen-Macaulay complexes $\G$ with $\F(\G) \subseteq \F(\D)$. Using standard techniques of linear programming we may express their inconsistency in terms of the exponents of the primary components of $I$. This approach was used before to study tetrahedral curves in \cite{GH}.
\smallskip

In Sections 2 and 3 we will use the above criteria to study the Cohen-Macaulayness of symbolic powers of the Stanley-Reisner ideal $I_\D$ of a simplicial complex $\D$. We will see that the Cohen-Macaulayness of $I_\D^{(2)}$ or of $I_\D^{(m)}$ for all $m \ge 1$ can be characterized completely in terms of $\D$ and that there are large classes of Stanley-Reisner ideals with Cohen-Macaulay symbolic powers.  
\smallskip

For every subset $V \subseteq [n]$ we denote by $\D_V$ the subcomplex of $\D$ whose facets are the facets of $\D$ with at least $|V|-1$ vertices in $V$. 
\medskip

\noindent{\bf Theorem 2.1.}
$I_\D^{(2)}$ is a Cohen-Macaulay ideal if and only if $\D$ is Cohen-Macaulay and $\D_V$ is Cohen-Macaulay for all subsets $V \subseteq [n]$ with $2 \le |V| \le \dim \D +1$.
\medskip

The condition of this theorem implies that the simplicial complex $\D$ is very compact in the sense that its vertices are almost directly connected to each other. In fact, we can show that the graph of the one-dimensional faces of $\D$ must have diameter $\le 2$. If $\D$ is a graph, we recover the result of \cite{MT} that  $I_\D^{(2)}$ is Cohen-Macaulay if and only if  the diameter of the graph is $\le 2$. Moreover, we also introduce a large class of simplicial complexes which generalizes matroid and shifted complexes and for which $I_\D^{(2)}$ is Cohen-Macaulay.
\medskip

In particular, using tools from linear programming we can show that the Cohen-Macaulayness of all symbolic powers characterizes matroid complexes.
\medskip

\noindent{\bf Theorem 3.5.}
$I_\D^{(m)}$ is Cohen-Macaulay for all $m \ge 1$ if and only if $\D$ is a matroid complex.
\medskip

This characterization is also proved independently by Varbaro \cite{Va}, who uses a completely different technique.
Theorem 3.5 adds a new algebraic feature to matroids, and we may hope that it could be used to obtain combinatorial information. As an immediate consequence we obtain the result of \cite{MT} that for  a graph $\D$, $I_\D^{(m)}$ is Cohen-Macaulay for all $m \ge 1$ if and only if every pair of disjoint edges of $\D$ is contained in a rectangle. Moreover, we can also easily deduce one of the main results of \cite{RTY} that for a flag complex $\D$, $I_\D^{(m)}$ is Cohen-Macaulay for all $m \ge 1$ if and only if the graph of the minimal nonfaces of $\D$ is a union of disjoint complete graphs.
\smallskip

It was showed in \cite{MT} and \cite{RTY} that if $\D$ is a graph or a flag complex and if $I_\D^{(t)}$ is Cohen-Macaulay for some $t \ge 3$, then $I_\D^{(m)}$ is Cohen-Macaulay for all $m \ge 1$. So one may ask the following general questions:

\begin{itemize}
\item Is $I_\D^{(m)}$ Cohen-Macaulay if $I_\D^{(m+1)}$ is Cohen-Macaulay?\par
\item Does there exist a number $t$ depending on $\dim \D$ such that if $I_\D^{(t)}$ is Cohen-Macaulay, 
then $I_\D^{(m)}$ is Cohen-Macaulay for all $m \ge 1$?
\end{itemize}

\noindent We don't know any definite answer to both questions. However, in Section 4 of this paper we can prove the following positive results on the preservation of Cohen-Macaulayness of symbolic powers.
\medskip

\noindent{\bf Theorem 4.3.}
$I_\D^{(m)}$ is Cohen-Macaulay if $I_\D^{(t)}$ is Cohen-Macaulay for some $t \ge (m-1)^2+1$.
\medskip

This result has the interesting consequence that $I_\D^{(2)}$ is Cohen-Macaulay if $I_\D^{(t)}$ is Cohen-Macaulay for some $t \ge 3$
or $I_\D^{(3)}$ is Cohen-Macaulay if $I_\D^{(t)}$ is Cohen-Macaulay for some $t \ge 5$. Note that we already know by \cite{HTT} that $I_\D$ is Cohen-Macaulay if $I_\D^{(t)}$ is Cohen-Macaulay for some $t \ge 2$.
\medskip

\noindent{\bf Theorem 4.5.}
Let $d = \dim \D$. If $I_\D^{(t)}$ is Cohen-Macaulay for some $t \ge (n-d)^{n+1}$, then 
$I_\D^{(m)}$ is Cohen-Macaulay for all $m \ge 1$.
\medskip

It remains to determine the smallest number $t_0$ such that if $I_\D^{(t)}$ is Cohen-Macaulay for some $t \ge t_0$, then $I_\D^{(m)}$ is Cohen-Macaulay for all $m \ge 1$. By \cite{MT} and \cite{RTY} we know that $t_0 = 3$ if $\dim \D = 1$ or if $\D$ is a flag complex. 
\smallskip

One may also raise similar questions on the Cohen-Macaulayness of ordinary powers of the Stanley-Reisner ideal $I_\D$.
Since $I_\D^m$ is Cohen-Macaulay if and only if $I_\D^{(m)} = I_\D^m$ and $I_\D^{(m)}$ is Cohen-Macaulay,
we have to study further the problem when $I_\D^{(m)} = I_\D^m$ in terms of $\D$.
The case $\dim \D = 1$ has been solved in \cite{MT}. We don't address this problem here because
it is of different nature than the Cohen-Macaulayness of $I_\D^{(m)}$ \cite{HHT}, \cite{HHTZ}.  
\smallskip

For unexplained terminology we refer the readers to the books \cite{BH}, \cite{Sch2} and \cite{Sta}.
\smallskip

Finally, the authors would like to thank the referee for suggesting Corollary \ref{shifted} and other corrections.

\section{Criteria for Cohen-Macaulay monomial ideals}

From now on let $I$ be a monomial ideal in the polynomial ring $S = k[x_1,...,x_n]$. 
Note that $S/I$ is an $\NN^n$-graded algebra.
For every degree $\a \in \ZZ^n$ we denote by $H_\mm^i(S/I)_\a$ the $\a$-component of the $i$-th local cohomology module $H_\mm^i(S/I)$ of $S/I$ with respect to the maximal homogeneous ideal $\mm$ of $S$.
Inspired of a result of Hochster in the squarefree case \cite[Theorem 4.1]{Ho}, Takayama found the following combinatorial formula for $\dim_kH_\mm^i(S/I)_\a$ \cite[Theorem 2.2]{Ta}.
\smallskip

For every $\a = (a_1,...,a_n) \in \ZZ^n$ we set $G_\a = \{i|\ a_i < 0\}$ and we denote by 
$\D_\a(I)$ the simplicial complex of all sets of the form $F \setminus G_\a$, where $F$ is a subset of $[n]$ containing $G_\a$ such that for every minimal generator $x^\b$ of $I$ there exists an index $i \not\in F$ such that $a_i < b_i$. Let $\D(I)$ denote the simplicial complex such that $\sqrt{I}$ is the Stanley-Reisner ideal of $\D(I)$. For simplicity we set $\D_\a = \D_\a(I)$ and $\D = \D(I)$.
\smallskip

For $j = 1,...,n$, let $\rho_j(I)$ denote the maximum of the $j$th coordinates of all vectors $\b \in \NN^n$ such that $x^\b$ is a minimal generator of $I$. 

\begin{Theorem}[Takayama's formula]\label{Takayama}
$$\dim_kH_\mm^i(S/I)_\a = 
\begin{cases}
\dim_k\tilde H_{i-|G_\a|-1}(\D_\a,k) & \text{\rm if }\ G_\a \in \D\ \text{\rm and}\\
&\ \ \  \ a_j  < \rho_j(I)\ \text{\rm for}\ j=1,...,n,\\
0 & \text{\rm else. }
\end{cases} $$
\end{Theorem}

It is known that $S/I$ is Cohen-Macaulay if and only if $H_\mm^i(S/I) = 0$ for $i < d$, where $d = \dim S/I$.
Therefore, we can derive from this formula criteria for the Cohen-Macaulayness of $I$.
The problem here is to find conditions by means of the primary decomposition of $I$.
The idea for that comes from \cite[Section 1]{MT}.
\smallskip

First, we have to describe the simplicial complexes $\D_\a$ in a more simple way.
For every subset $F$ of $[n]$ let $S_F = S[x_i^{-1}|\ i \in F]$.

\begin{Lemma}\label{face} 
$\D_\a$ is the simplicial complex of all sets of the form $F \setminus G_\a$, where $F$ is a subset of $[n]$ containing $G_\a$ such that $x^\a \not\in IS_F$.
\end{Lemma}

\begin{proof}
We have $a_i < b_i$ for some $i \not\in F$ iff $x^\a$ is not divided by $x^\b$ in $S_{F}$. 
This condition is satisfied for every minimal generator $x^\b$ of $I$ iff $x^\a \not\in IS_F$.
\end{proof}

This lemma can be also proved by looking at the $\a$th multigraded component of the {\v C}ech complex of $S/I$. \smallskip

Using the above characterization of $\D_\a$ we can easily show that $\D_\a$ is a subcomplex of $\D$.
In fact, $\D$ is the simplicial complex of all subsets $F \subseteq [n]$ such that $\prod_{i \in F}x_i \not\in \sqrt{I}$. But this condition means $IS_F \neq S_F$.
If $G \in \D_\a$, then $x^\a \not\in IS_G$, which implies $IS_G \neq S_G$, hence $G \in \D$. This shows that $\D_\a \subseteq \D$.

\begin{Example}\label{0}
{\rm $\D_\0 = \D$ because for all faces $F$ of $\D$ we have $x^\0 = 1 \not\in IS_F$.}
\end{Example}

For every subset $F$ of $[n]$ let $P_F$ denote the prime ideal of $S$ generated by the variables $x_i$, $i \not\in F$. 
Then the minimal primes of $I$ are the ideals $P_F$, $F \in \F(\D)$.
Let $I_F$ denote the $P_F$-primary component of $I$. If $I$ has no embedded components, we have
$$I = \bigcap_{F \in \F(\D)}I_F.$$
Using this primary decomposition of $I$ we obtain the following formula for the dimension of $\D_\a$.

\begin{Lemma}\label{dim}
Assume that $I$ is unmixed. Then $\D_\a(I)$ is pure and
$$\dim \D_\a = \dim \D - |G_\a|.$$
\end{Lemma}

\begin{proof} 
The assumption means that $I$ has no embedded components and $\D$ is pure. Let $H$ be an arbitrary facet of $\D_\a$.  By Lemma \ref{face}, $x^\a \not\in IS_{H \cup G_\a}$. We have
$$IS_{H \cup G_\a}  =  \bigcap_{F \in \F(\D)}I_FS_{H \cup G_\a}= \bigcap_{F \in \F(\D),\, H \cup G_\a \subseteq F}I_FS_{H \cup G_\a}$$
because $P_FS_{H \cup G_\a} = S_{H \cup G_\a}$ if $H \cup G_\a \not\subseteq F$. Therefore, there exists $F \in \F(\D)$ with $H \cup G_\a \subseteq F$ such that $x^\a \not\in I_FS_{H \cup G_\a}$. Since $I_FS_{H \cup G_\a} \cap S_F = IS_F$, this implies $x^\a \not\in IS_F$, hence $F \setminus G_\a \in \D_\a$ by Lemma \ref{face}. So we must have $H = F \setminus G_\a$. Thus, 
$$\dim H = |F| - |G_\a| - 1 = \dim \D - |G_\a|.$$
This shows that $\D_\a$ is pure and $\dim \D_\a = \dim \D - |G_\a|.$
\end{proof}

If $\a \in \NN^n$, then $G_\a = \emptyset$. Hence Lemma \ref{dim} implies $\F(\D_\a) \subseteq \F(\D)$. 
We can easily check which facet of $\F(\D)$ belongs to $\F(\D_\a)$ and we can determine all points $\a \in \NN^n$ such that $\D_\a$ equals to a given subcomplex $\G$ of $\D$ with $\F(\G) \subseteq \F(\D)$. 
For that we introduce the set of lattice points
$$L_\G(I) := \big\{\a \in \NN^n|\  x^\a \in \bigcap_{F \in \F(\D)\setminus\F(\G)}I_F \setminus \bigcup_{G\in \F(\G)}I_G\big\}.$$

\begin{Lemma}\label{division}
Assume that $I$ is unmixed. For $\a \in \NN^n$ we have \par
{\rm (i) } $\F(\D_\a) = \big\{F \in \F(\D)|\ x^\a \not\in I_F\big\}$,\par
{\rm (ii) } $\D_\a = \G$ if and only if $\a \in L_\G(I)$.
\end{Lemma}

\begin{proof} 
For $F, G \in \F(\D)$ we have $I_GS_F = S_F$ if $G \neq F$. Therefore,
$$IS_F \cap S = \bigcap_{G \in \F(\D)} I_GS_F \cap S = I_FS_F \cap S = I_F.$$
From this it follows that $x^\a \in IS_F$ iff $x^\a \in I_F$. By  Lemma \ref{face}, $F \in \F(\D_\a)$ iff $x^\a \not\in I_F$,
which immediately yields the assertions. 
\end{proof}

With regard to Lemma \ref{division} we may consider the following two criteria for the Cohen-Macaulayness of $I$ as by means of the primary decomposition of $I$. 
\smallskip

For any face $F$ of a simplicial complex $\G$ we denote by $\lk_\G F$ the subcomplex of all faces $G \in \G$ such that $F \cap G = \emptyset$ and $F \cup G \in \G$. We call $\G$ a {\it Cohen-Macaulay complex} (over $k$) if $\tilde H_j(\lk_\G F,k) = 0$ for all $F \in \G$, $j < \dim \lk_\G F$.

\begin{Theorem}\label{criteria}
Assume that $I$ is an unmixed monomial ideal. Then the following conditions are equivalent:\par
{\rm (i)} $I$ is a Cohen-Macaulay ideal,\par
{\rm (ii)} $\D_\a$ is a Cohen-Macaulay complex for all $\a \in \NN^n$,\par
{\rm (iii)} $L_\G(I) = \emptyset$ for every non-Cohen-Macaulay complex $\G$ with $\F(\G) \subseteq \F(\D)$.
\end{Theorem}

\begin{proof}
(i) $\Rightarrow$ (ii): 
Let $F \in \D_\a$ be arbitrary. We will first represent $\lk_{\D_\a}F$ in a suitable form in order to apply Takayama's formula.
Let $G \in \D_\a$ such that $F \cap G = \emptyset$. By Lemma \ref{face}, $F \cup G \in \D_\a$ iff $x^\a \not\in IS_{F\cup G}$.
Let $\b \in \ZZ^n$ such that $b_i = -1$ for $i \in F$ and $b_i = a_i$ for $i \not\in F$. Then $F = G_\b$, and $x^\a \not\in IS_{F\cup G}$ iff $x^\b \not\in IS_{F \cup G}$. By Lemma \ref{face}, $G \in \D_\b$ iff $x^\b \not\in IS_{F \cup G}$. Therefore, $F \cup G \in \D_\a$ iff $G \in \D_\b$. So we obtain $\lk_{\D_\a}F = \D_\b$. By the proof of \cite[Theorem 1]{Ta}, $\tilde H_i(\D_\b,k) = 0$ for all $i$ if there is a component $b_j \ge \rho_j(I)$. Therefore, we may assume that
 $b_j < \rho_j(I)$ for all $j$. By Theorem \ref{Takayama} the Cohen-Macaulayness of $I$ implies 
$$\tilde H_{i-|G_\b|-1}(\D_\b,k) = 0\ \text{for}\ i < d,$$
where $d = \dim S/I$. By Lemma \ref{dim}, 
$$\dim \D_\b = \dim \D - |G_\b| = d - |G_\b| - 1.$$
Therefore, the above formula can be rewritten as
$$\tilde H_j(\D_\b,k) = 0\ \text{for}\ j < \dim \D_\b.$$
So we can conclude that $\tilde H_j(\lk_{\D_\a}F,k) = 0$ for $j < \dim \lk_{\D_\a}F.$\par

(ii) $\Rightarrow$ (iii): 
By Lemma \ref{division}(ii), $\D_\a(I) = \G$ for all $\a \in L_\G(I)$. Therefore, $L_\G(I) = \emptyset$ if $\G$ is not Cohen-Macaulay.\par

(iii) $\Rightarrow$ (i): 
By Theorem \ref{Takayama} we only need to show that 
$\tilde H_{i-|G_\a|-1}(\D_\a,k) = 0$ 
for all $\a \in \ZZ^n$ with $G_\a \in \D$, $i < d$. As we have seen above, this formula can be rewritten as 
$\tilde H_j(\D_\a,k) = 0$ for $j < \dim \D_\a.$
We may assume that $\D_\a \neq \emptyset$. By Lemma \ref{face}, there is a set $G \supseteq G_\a$ such that $x^\a \not\in IS_G$. From this it follows that $x^\a \not\in IS_{G_\a}$.
Let $\b \in \NN^n$ with $b_i = a_i$ if $a_i \ge 0$ and $b_i = 0$ else. Then $x^\b \not\in IS_F$ iff $x^\a \not\in IS_F$, $F \supseteq G_\a$. So
$x^\b \not\in IS_{G_\a}$. Let $\G = \D_\b$. 
By Lemma \ref{face}, $G_\a \in \G$ and
$$\D_\a = \{F \setminus G_\a|\ F \supseteq G_\a, x^\b \not\in IS_F\} = \{F \setminus G_\a|\  F \supseteq G_\a, F \in \G\}
= \lk_{\G}G_\a.$$ 
By Lemma \ref{division}, $\F(\G) \subseteq \F(\D)$ and $\b \in L_\G(I)$. Therefore, (iii) implies that $\G$ is Cohen-Macaulay. Hence
$\tilde H_j(\D_\a,k) = 0$ for $j < \dim \D_\a.$
\end{proof}

\begin{Remark}\label{restriction}
{\rm The above proof also shows that we may replace Theorem \ref{criteria}(ii) by the condition that  $\D_\a$ is Cohen-Macaulay for $\a \in \NN^n$ with $a_j < \rho_j(I)$, $j = 1,...,n$. This restriction is very useful in computing examples.}
\end{Remark}

If $I$ is a squarefree monomial ideal, $\rho_j(I) = 1$ for all $j$, hence there is only a point $\a \in \NN^n$ with $a_j < 1$ for all $j$, which is $\0$. But $\D_\0 = \D$. Therefore, Theorem \ref{criteria}(ii) implies the well-known result that $I$ is Cohen-Macaulay if and only if $\D$ is Cohen-Macaulay \cite{Re}. If $I$ is an arbitrary monomial ideal, Theorem \ref{criteria}(ii) implies that $\D$ is Cohen-Macaulay if $I$ is Cohen-Macaulay. From we immediately obtain the result that $\sqrt{I}$ is Cohen-Macaulay \cite[Theorem 2.6(i)]{HTT}.
\smallskip

If $I$ is the intersection of prime ideal powers, we can interpret Theorem \ref{criteria}(iii) in terms of Diophantine linear inequalities. In fact, if $I_F = P_F^{m_F}$ for some positive integer $m_F$, we have $x^a \in I_F$ if and only if $\sum_{i \not\in F} a_i \ge m_F$. Hence we can translate the condition
$$x^\a \in \bigcap_{F \in \F(\D) \setminus \F(\G)}I_F \setminus \bigcup_{G \in \F(\G)}I_G$$
as a system of linear inequalities:
\begin{align*}
\sum_{i \not\in F}a_i \ge m_F &\  \big(F \in \F(\D) \setminus \F(\G)\big),\\
\sum_{i \not\in G}a_i < m_G &\ \big(G \in \F(\G)\big).
\end{align*}
The condition $L_\G(I) = \emptyset$ means that this system of linear inequalities has no solution $\a \in \NN^n$. Thus, $I$ is Cohen-Macaulay if and only if this system is inconsistent in $\NN^n$ for all non-Cohen-Macaulay subcomplexes $\G$ of $\D$ with $\F(\G) \subseteq \F(\D)$.\smallskip

In particular, if $\dim S/I = 2$, we may identify $\D$ with the graph of its edges. In this case, the non-Cohen-Macaulay subcomplexes are the unconnected subgraphs so that we can easily write down the corresponding systems of linear inequalities. As an example we consider the following class of monomial ideals for which it took several efforts \cite{Schw}, \cite{MN} until one knows which of them is Cohen-Macaulay \cite{Fr}, \cite{GH}. 

\begin{Example}[Tetrahedral curves] {\rm
Let 
$$I = (x_1,x_2)^{m_1}\cap (x_1,x_3)^{m_2} \cap (x_1,x_4)^{m_3} \cap (x_2,x_3)^{m_4} \cap (x_2,x_4)^{m_5} \cap (x_3,x_4)^{m_6},$$
where $m_1,...,m_6$ are arbitrary positive integers. Then $\D$ is the complete graph $K_4$.   
This graph has three unconnected subgraphs which correspond to the pairs of disjoint edges:
$\left\{\{1,2\},\{3,4\}\right\}, \left\{\{1,3\},\{2,4\}\right\}, \left\{\{1,4\},\{2,3\}\right\}.$
Let $\G$ be the complex of the subgraph $\left\{\{1,2\},\{3,4\}\right\}$. Then $L_\G(I)$ is the set of all points $\a  \in \NN^4$ which satisfies the inequalities
$$\left\{\begin{array}{ll}
a_2 + a_4 \ge m_2,&   a_2 + a_3 \ge m_3,\ \ a_1 + a_4 \ge m_4,\ \  a_1 + a_3 \ge m_5,\\
a_3+a_4 < m_1,&  a_1+a_2 < m_6.
\end{array}\right.$$
For the complexes of the subgraphs $\left\{\{1,3\},\{2,4\}\right\}$, $\left\{\{1,4\},\{2,3\}\right\}$ we have two similar systems of inequalities.
By Theorem \ref{criteria}(iii), $I$ is Cohen-Macaulay iff the three systems of inequalities have no solutions in $\NN^4$. Using standard techniques of integer programming one can easily solve these systems of inequalities and obtain a Cohen-Macaulay criterion for $I$ in terms of the exponents $m_1,...,m_6$ (see \cite{GH} for details).}
\end{Example}

Recently, Herzog, Takayama and Terai \cite[Theorem 3.2]{HTT} proved that all unmixed monomial ideals with radical $I_\D$ are Cohen-Macaulay if and only if $\D$ has no non-Cohen-Macaulay subcomplex $\G$ with $\F(\G) \subseteq \F(\D)$. But that is just an immediate consequence of Theorem \ref{criteria}(iii). 
In addition, we can use the same condition on $\D$ to characterize the Cohen-Macaulayness of all intersections of prime ideal powers with radical $I_\D$.

\begin{Corollary}
Let $\D$ be a pure simplicial complex. The ideal $I = \cap_{F \in \F(\D)}P_F^{m_F}$
is Cohen-Macaulay for all exponents $m_F \ge 1$ (or $m_F \gg 0$) if and only if $\D$  has no non-Cohen-Macaulay subcomplex $\G$ with $\F(\G) \subseteq \F(\D)$. 
\end{Corollary}

\begin{proof} 
It suffices to show the necessary part. Assume that $\D$ has a non-Cohen-Macaulay subcomplex $\F$ with $\F(\G) \subseteq \F(\D)$. Given any point $\a \in \NN^n$ with all $a_i >0$ we choose
\begin{align*}
m_F & = \sum_{i \not\in F}a_i \ \  \big(F \in \F(\D) \setminus \F(\G)\big),\\
m_G & = \sum_{i \not\in G}a_i +1\ \ \big(G \in \F(\G)\big).
\end{align*}
As mentioned above, this implies $L_\G(I) \neq \emptyset$. Hence $I$ is not Cohen-Macaulay by Theorem \ref{criteria}(iii). 
\end{proof}

Note that $\D$  has no non-Cohen-Macaulay subcomplex $\G$ with $\F(\G) \subseteq \F(\D)$ if and only if after a suitable permutation, $\F(\D) = \{F_1,...,F_r\}$ with $F_i = \{1,...,i-1,i+1,...,d+1\}$, $i = 1,...,r$, or $F_i = \{1,...,d,d+i\}$, $i = 1,...,r$ \cite[Theorem 3.2]{HTT}.

\section{Cohen-Macaulayness of the second symbolic power}

Let $\D$ be an arbitrary simplicial complex on the vertex set $[n]$. One calls 
$$I_\D = \bigcap_{F \in \F(\D)}P_F,$$
the Stanley-Reisner ideal and $k[\D] = S/I_\D$ the face ring of $\D$. 
\smallskip

We will use Theorem \ref{criteria} to study the Cohen-Macaulayness of the symbolic powers of $I_\D$. 
For every integer $m \ge 1$, the $m$-th {\it symbolic power} of $I_\D$ is the ideal
$$I_\D^{(m)} =  \bigcap_{F \in \F(\D)}P_F^m.$$
Obviously, $\D(I_\D^{(m)}) = \D$. Since we study the Cohen-Macaulayness of $I_\D^{(m)}$ we may assume that $\D$ is pure, which is equivalent to say that $I_\D$ is  unmixed.
\smallskip

For every subset $V \subseteq [n]$ we denote by $\D_V$ the subcomplex of $\D$ whose facets are the facets of $\D$ with at least $|V|-1$ vertices in $V$. 

\begin{Theorem}\label{second}
$I_\D^{(2)}$ is Cohen-Macaulay if and only if $\D$ is Cohen-Macaulay and $\D_V$ is Cohen-Macaulay for all subsets $V \subseteq [n]$ with $2 \le |V| \le \dim \D +1$.
\end{Theorem}

\begin{proof}
For $I =  I_\D^{(2)}$ we have $\rho_j(I) = 2$ for all $j= 1,..,n$. Hence $\{0,1\}^n$ is the set of all $\a \in \NN^n$ with $a_j < \rho_j(I) = 2$, $j = 1,..,n$. 
By Remark \ref{restriction}, $I_\D^{(2)}$ is Cohen-Macaulay iff $\D_\a$ is Cohen-Macaulay for all $\a \in \{0,1\}^n$. \par
If $\a = 0$, $\D_\0 = \D$ by Example \ref{0}. If $\a = \e_1,...,\e_n$, the unit vectors of $\NN^n$, we have 
$x^{\e_i} = x_i \not\in P_F^2$ for all $F \in \F(\D)$, which by Lemma \ref{division}(i) implies $\D_\a = \D$.
If $\a \neq 0,\e_1,...,\e_n$, let $V = \{i \in [n]|\ a_i = 1\}$. Then $|V| \ge 2$ and $\D_\a = \D_V$. 
In fact, for any subset $F$ of $[n]$, $F$ is a facet of $\D_\a$ iff $x^\a \not\in P_F^2$ iff $\sum_{i \not\in F}a_i < 2$ 
iff $|V \setminus F| < 2$ iff $|F \cap V| \ge |V|-1$. \par
It remains to show that $\D_V$ is Cohen-Macaulay if $|V| \ge \dim \D+2$. 
If $|V| = \dim \D +2$, then $\D_V$ is a union of facets of a simplex. In this case, $I_{\D_V}$ is a principal ideal. Hence $\D_V$ is Cohen-Macaulay. 
If $|V| \ge \dim \D + 3$, then $\D_V = \emptyset$ because no facet of $\D$ can have more than $\dim \D+1$ vertices. 
\end{proof}

Theorem \ref{second} puts strong constraints on simplicial complexes $\D$ for which $I_\D^{(2)}$ is Cohen-Macaulay. 
We shall see later in Corollary \ref{down to 2} that $I_\D^{(2)}$ is Cohen-Macaulay if $I_\D^{(m)}$ is Cohen-Macaulay for some $m \ge 3$. 
\smallskip

Recall that for a graph $\G$, the distance between two vertices of $\G$ is the minimal length of paths from one vertex to the other vertex. This length is infinite if there is no paths connecting them. The maximal distance between two vertices of $\G$ is called the {\it diameter} of $\G$ and denoted by $\diam(\G)$.  

\begin{Corollary}\label{diameter}
Let $\D$ be a simplicial complex such that $I_\D^{(2)}$ is a Cohen-Macaulay ideal. 
Let $\G$ be the graph of the one-dimensional faces of $\D$. Then $\diam(\G) \le 2$.
\end{Corollary}

\begin{proof}
Let $i \neq j$ be two arbitrary vertices of $\G$ and put $V = \{i,j\}$. Then the faces of $\D_V$ are the faces of $\G$ which contain $i$ or $j$. By Theorem \ref{second}, $\D_V$ is connected. Therefore, there are a face containing $i$ and a face containing $j$ which meet each other. This implies that $\G$ has  an edge containing $i$ and an edge containing $j$ which share a common vertex. Hence
the distance between $i$ and $j$ is $\le 2$. 
\end{proof}

The converse of Corollary \ref{diameter} holds in the case $\dim \D = 1$. 

\begin{Corollary}\label{graph} {\rm \cite[Theorem 2.3]{MT}}
Let $\D$ be a graph. Then $I_\D^{(2)}$ is a Cohen-Macaulay ideal if and only if $\diam(\D) \le 2$.
\end{Corollary}

\begin{proof}
It is known that a graph is Cohen-Macaulay iff it is connected.
Therefore, it suffices to show that $\D_V$ is connected for all subsets $V \subseteq [n]$ with $|V| = 2$ iff $\diam(\D) \le 2$.
Assume that $V = \{i,j\}$. Then the edges of $\D_V$ are the edges of $\D$ which contain $i$ or $j$. Therefore, $\D_V$ is connected iff the distance between $i$ and $j$ is $\le 2$. Since $i,j$ can be chosen arbitrarily, this means $\diam(\D) \le 2$.
\end{proof}

Munkres \cite{Mu} showed that the Cohen-Macaulayness of $I_\D$ depends only on the geometric realization of $\D$.
In other words, the Cohen-Macaulayness of $I_\D$ is a topological property of $\D$.
From Corollary \ref{graph} we can easily see that the Cohen-Macaulayness of $I_\D^{(2)}$ is not a topological property of $\D$.

\begin{Example}
{\rm Let $\D$ be a path of length $r$. Then $\diam(\D) = r$. Hence $I_\D^{(2)}$ is Cohen-Macaulay if $r = 1,2$ and not Cohen-Macaulay if $r \ge 3$, though
any path is topologically a line. Since the barycentric subdivision of a path of length 2 is a path of length 4, this also shows that  the Cohen-Macaulayness of $I_\D^{(2)}$ doesn't pass to the barycentric subdivision of $\D$.}
\end{Example}

For higher dimensional simplicial complexes we couldn't get a similar result as Corollary \ref{graph} because we don't know how to check the Cohen-Macaulayness of subcomplexes. This can be done only in special cases. 
\smallskip

We call a pure simplicial complex $\D$ a {\it tight complex} if there is a labelling of the vertices such that  for every pair of facets $G_1, G_2$ and vertices $i \in G_1 \setminus G_2, j \in G_2 \setminus G_1$ with $i < j$ there is a vertex $j' \in G_1 \setminus G_2$ such that $(G_2  \setminus \{j\}) \cup \{j'\}$ is a facet.  Obviously, this class of complexes contains all matroid complexes. \smallskip

Recall that a {\it matroid complex} is a collection of subsets of a finite set, called {\it independent sets}, with the following properties:
\smallskip

	(1) The empty set is independent.\par
   (2) Every subset of an independent set is independent. \par
   (3) If $F$ and $G$ are two independent sets and $F$  has more elements than $G$, then there exists an element in $F$ which is not in $G$ that when added to $G$ still gives an independent set. 
\smallskip

Examples of matroid complexes are abundant such as collections of linearly independent subsets of finite sets of elements in a vector space.
Note that there are tight complexes of any dimension which are not matroid complexes such as the complex generated by all subsets of $n-2$ elements of $[n-1]$ and the set $\{3,...,n\}$, $n \ge 4$.

\begin{Theorem}\label{tight}
Let $\D$ be a tight complex. Then $I_\D^{(2)}$ is Cohen-Macaulay.
\end{Theorem}

\begin{proof}
If $n = 2$, the assertion is trivial. Assume that $n > 2$.
By Theorem \ref{criteria}(iii) we only need to show that $L_\G(I_\D^{(2)}) = \emptyset$ for all non-Cohen-Macaulay subcomplexes $\G$ of $\D$ with $\F(\G) \subseteq \F(\D)$. Without restriction we may assume that $n \in \G$. \par

Let $\D_1$  and $\G_1$ be the subcomplexes of $\D$ and $\G$ generating by the facets not containing $n$, respectively.
Then  $\D_1$ is a tight complex on $[n-1]$ and $\G_1$ is a subcomplex of $\D_1$ with $\F(\G_1) \subseteq \F(\D_1)$. 
Since $\F(\G_1) \subseteq \F(\G)$ and $\F(\D_1) \setminus \F(\G_1) \subseteq \F(\D) \setminus \F(\G)$, $L_\G(I_\D^{(2)}) \subseteq L_{\G_1}(I_{\D_1}^{(2)})$.
By induction we may assume that $I_{\D_1}^{(2)}$ is Cohen-Macaulay. 
If $\G_1$ is not Cohen-Macaulay, $L_{\G_1}(I_{\D_1}^{(2)}) = \emptyset$ by Theorem \ref{criteria}(iii) and hence $L_\G(I_\D^{(2)}) = \emptyset$.
So we may assume that $\G_1$ is Cohen-Macaulay. \par

Let $\D_2 = \lk_\D \{n\}$. Then $\D_2$ is a tight complex on $[n-1]$. By induction we may assume that $I_{\D_2}^{(2)}$ is Cohen-Macaulay.  Let $\D_2^*$ be the subcomplex of $\D$ generating by the facets containing $n$. Then $I_{\D_2}$ and $I_{\D_2^*}$ lie in different polynomial rings but have the same (minimal) monomial generators. Therefore, $I_{\D_2^*}^{(2)}$ is Cohen-Macaulay. \par

Let $\G_2^*$ be the subcomplex of $\G$ generating by the facets containing $n$. Then $\G_2^*$ is a subcomplex of $\D_2^*$ with $\F(\G_2^*) \subseteq \F(\D_2^*)$. 
If $\G_2^*$ is not Cohen-Macaulay,  $L_{\G_2^*}(I_{\D_2^*}^{(2)}) = \emptyset$ by Theorem \ref{criteria}(iii). On the other hand, it is easy to see that  $L_\G(I_\D^{(2)}) \subseteq L_{\G_2^*}(I_{\D_2^*}^{(2)})$.
Therefore, $L_\G(I_\D^{(2)}) = \emptyset$. So we may assume that $\G_2^*$ is Cohen-Macaulay. \par 

Let $\G_2 = \lk_\G\{n\}$. Since $\G_2^*$ is a cone over $\G_2 $, $\G_2 $ is Cohen-Macaulay.
Note that $\G_1 \cap \G_2^* \subseteq \G_2$ and $\G_1 \cup \G_2^* = \G$. 
If $\G_1 \cap \G_2^* = \G_2$, there is an exact sequence
$$0 \to k[\G] \to k[\G_1]\oplus k[\G_2^*] \to k[\G_2] \to 0.$$
Since $k[\G_1], k[\G_2^*]$ and $k[\G_2]$ are Cohen-Macaulay
with $\dim k[\G_1] = \dim k[\G_2^*] = \dim k[\G_2] +1$, 
we can conclude that $k[\G]$ is Cohen-Macaulay, which  contradicts the assumption that $\G$ is not Cohen-Macaulay.
So $\G_1 \cap \G_2^*$ is properly contained in  $\G_2$.  \par

This means that there exists a facet $G \in \F(\G)$ containing $n$ such that $G \setminus \{n\}$ is not contained in any facet of $\F(\G)$ not containing $n$. Moreover, there also exists a facet of $\F(\G)$ not containing $n$ because otherwise $\G = \G_2^*$ were Cohen-Macaulay.  By the definition of tight complexes we can see that these properties hold for any vertex. \par

Assume for the contrary that  $L_\G(I_\D^{(2)}) \neq \emptyset$ and choose $\a \in L_\G(I_\D^{(2)})$ arbitrary. By the proof of Theorem 2.1, $\a \in \{0,1\}^n$ with $|\{i \in [n]|\ a_i = 1\}| \le \dim \D+1$. Since  $n > \dim \D+1$, there is at least a vertex $j$ with $a_j = 0$. Let $j = \max\{i \in [n]|\ a_i = 0\}$. \par

Choose  a facet $G_1 \in \F(\G)$ not containing $j$ and a facet $G_2 \in \F(\G)$ containing $j$ such that $G_2 \setminus \{j\}$ is not contained in any facet of $\F(\G)$ not containing $j$.  If there is a vertex $i \in G_1 \setminus G_2$ such that $i < j$, there is a vertex $j' \in G_1 \setminus G_2$ such that $F = (G_2 \setminus \{j\}) \cup \{j'\}$ is a facet of $\D$. By the choice of $G_2$, $F \not\in \F(\G)$. So we have
$\sum_{i \not\in F}a_i \ge 2 \  \text{and}\ \sum_{i \not\in G_2}a_i < 2.$
From this it follows that $a_j > a_{j'}$, which is a contradiction because $a_j = 0$ and $a_{j'} \ge 0$. Thus, $i > j$ and hence $a_i = 1$ for every vertex $i \in G_1 \setminus G_2$. Since $j \not\in G_1$ and $G_2  \setminus \{j\} \not\subseteq  G_1$, $|G_1 \cap G_2| \le |G_2| -2 = |G_1|-2$. Thus, $G_1 \setminus G_2$ contains at least two vertices, say $i$ and $i'$. Since $a_i = a_{i'} = 1$, we get
$\sum_{t \not\in G_2}a_t \ge a_i + a_{i'} = 2$, a contradiction.
So we have proved that $L_\G(I_\D^{(2)}) = \emptyset$.
\end{proof}

The converse of Theorem \ref{tight} is not true. 

\begin{Example}
{\rm Let $\D$ be the graph of a 5-cycle:

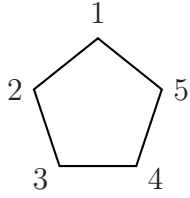
\begin{center}
\psset{unit=1.7cm}
\begin{pspicture}(-0.7,-0.8)(0.7,0.85)
\pspolygon(0,0.5)(-0.5,0.1)(-0.3,-0.5)(0.3,-0.5)(0.5,0.1)
\rput(0,0.7){1}
 \rput(-0.65,0.1){2}
 \rput(-0.45,-0.6){3}
 \rput(0.45,-0.6){4}
 \rput(0.65,0.1){5}
\end{pspicture}
\end{center}

\noindent Then $\diam(\D) = 2$ so that $I_\D^{(2)}$ is Cohen-Macaulay by Corollary \ref{graph}.
For any labelling of the vertices of $\D$ we consider an arbitrary pair of disjoint edges $\{i,i'\}$ and $\{j,j'\}$. 
Without restriction we may assume that $\{i,j\}$ is an edge of $\D$. Then $i'$ and $j'$ is not connected to the
edges $\{j,j'\}$ and $\{i,i'\}$ by any edge, respectively. Hence $\D$ is not a tight complex.}
\end{Example}

The proof of Theorem \ref{tight} is remarkable in the sense that it gives a method to pass the difficult test on all non-Cohen-Macaulay subcomplexes of $\D$ to the unions of two facets. We will use it again in the proof of Theorem \ref{all matroids}.
\smallskip

A simplicial complex $\D$ on the vertex set $[n]$ is called a {\it shifted complex} if there is a labelling of the vertices such that for every face $F \in \D$ and every vertex $i \in F$, $(F \setminus \{i\}) \cup \{j\} \in \D$ for all $j < i$ \cite{Ka}. Obviously, shifted complexes are tight.

\begin{Corollary} \label{shifted}
Let $\D$ be a pure shifted complex. Then $I_\D^{(2)}$ is Cohen-Macaulay.
\end{Corollary}

We now present an operation for the construction of new simplicial complexes such that the second symbolic power of their Stanley-Reisner ideals are Cohen-Macaulay. Given two simplicial complexes $\D$ and $\G$ on disjoint vertex sets, one calls the simplicial complex 
$$\D*\G = \{F \cup G|\ F \in \D, G \in \G\}$$
the {\it join} of $\D$ and $\G$. 

\begin{Theorem}
Let $\D$ and $\G$ be simplicial complexes such that $I_{\D}^{(2)}$ and $I_{\G}^{(2)}$ are Cohen-Macaulay.
Then $I_{\D*\G}^{(2)}$ is Cohen-Macaulay.
\end{Theorem}

\begin{proof}
Let $\D$ and $\G$ be complexes on the vertex sets $[n]$ and $\{n+1,...,n+m\}$, respectively.
Let $d = \dim \D$ and $e = \dim \G$. Then $\dim \D*\G = d + e + 1$.
By Theorem \ref{second} we have to show that $(\D*\G)_U$ is Cohen-Macaulay for all 
$U \subseteq [m+n]$ with $2 \le |U| \le d+e+2$. \par
Set $V = U \cap [n]$ and $W = U \cap \{n+1,...,n+m\}$. Let
\begin{eqnarray*}
\star_\D\! V & = & \{F \in \D|\ F \cup V \in \D\},\\ 
\star_\G\! W & = & \{G \in \G|\ G \cup W \in \G\}.
\end{eqnarray*}
It is easy to see that 
$$(\D*\G)_U = (\D_V*\star_\G\! W) \cup (\star_\D\! V*\G_W).$$
By Theorem \ref{second}, $\D_V$ and $\G_W$ are Cohen-Macaulay. By \cite[Chapter III, Proof of Corollary 9.2]{Sta}, $\star_\D\! V$ and $\star_\G W$ are Cohen-Macaulay. Therefore, $\D_V*\star_\G\! W$ and $\star_\D\! V*\G_W$ are Cohen-Macaulay complexes of dimension $d+e+1$ \cite[Exercise 5.1.21]{BH}.
Since $\star_\D\! V \subseteq \D_V$ and $\star_\G\! W \subseteq \G_W$,
$$ (\D_V*\star_\G\! W) \cap (\star_\D\! V*\G_W) = \star_\D\! V * \star_\G\! W.$$
Therefore,  $(\D_V*\star_\G\! W) \cap (\star_\D\! V*\G_W)$ is a Cohen-Macaulay complex of dimension $d+e-1$. Now, from the exact sequence
$$0 \to k[(\D*\G)_U] \to k[\D_V*\star_\G\! W]\oplus k[\star_\D\! V*\G_W] \to k[(\D_V*\star_\G\! W) \cap (\star_\D\! V*\G_W)] \to 0$$
we can conclude $k[(\D*\G)_U]$ is Cohen-Macaulay. 
\end{proof}

It is well known that the Cohen-Macaulayness of $I_\D$ depends on the characteristic of the base field \cite{Re}. By Theorem \ref{second} we may expect  that the Cohen-Macaulayness of $I_\D^{(2)}$ also depends on the characteristic of the base field. However we have been unable to settle this problem. The triangulation of the projective plane does not provide an example for that.

\begin{Example}
{\rm
Let $\D$ be the triangulation of the projective plane with the facets
\begin{align*} 
\{1,2,3\},\{1,2,6\},\{1,3,5\},\{1,4,5\},\{1,4,6\},\{2,3,4\},&\\
\{2,4,5\},\{2,5,6\},\{3,4,6\},\{3,5,6\}.&
\end{align*}

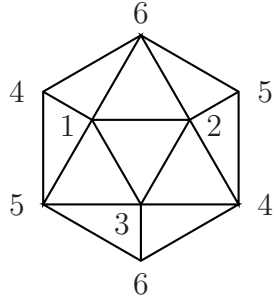
\begin{center}
\psset{unit=1.5cm}
\begin{pspicture}(-2,-1.3)(2,1.3)
\pspolygon(0,1)(-0.866,0.5)(-0.866,-0.5)
\pspolygon(0,1)(0.866,0.5)(0.866,-0.5)
 \pspolygon(0,-1)(-0.866,-0.5)(0.866,-0.5)
 \pspolygon(0.433,0.25)(0,-0.5)(-0.433,0.25)
 \psline(0.866,0.5)(0.433,0.25)
 \psline(-0.866,0.5)(-0.433,0.25)
 \psline(0,-0.5)(0,-1)
 \rput(-0.65,0.2){1}
 \rput(0.65,0.2){2}
 \rput(-0.17,-0.67){3}
 \rput(0,1.2){6}
 \rput(1.1,0.5){5}
 \rput(1.1,-0.5){4}
 \rput(0,-1.2){6}
 \rput(-1.1,-0.5){5}
 \rput(-1.1,0.5){4}
\end{pspicture}
\end{center}

\noindent Since all vertices of $\D$ are connected by one-dimensional faces, $\diam(G) = 1$. 
On the other hand, for $V = \{4,5,6\}$, $\D_V$ is the simplicial complex with the facets $\{1,4,5\},\{1,4,6\},\{2,5,6\},\{2,4,5\},\{3,4,6\},\{3,5,6\}.$ 
Since the geometric realization of $\D_V$ can be contracted to a cycle, $\D_V$ is not Cohen-Macaulay. Therefore,
$I_\D^{(2)}$ is not Cohen-Macaulay by Theorem \ref{second}. }
\end{Example}

\section{Cohen-Macaulayness of all symbolic powers}

In the following we shall use Theorem \ref{criteria}(iii) to study the Cohen-Macaulayness of all symbolic powers of Stanley-Reisner ideals.
For that we shall need the following characterization of strict homogeneous inequalities.

\begin{Lemma}\label{Motzkin}
Let $A$ and $B$ be matrices having the same number of columns. Then there exists a column vector $\x$ such that $A\x < \0$ and $B\x \ge \0$ if and only if there are no row vectors $\y, \z \ge \0$ such that $\y A + \z B = \0$ and $\y \neq \0$.
\end{Lemma}

\begin{proof}
Consider the general system $A\x < \b$ and $B\x \ge \c$, where $\b$ and $\c$ are given column vectors.
Motzkin's transposition theorem (see e.g. \cite[Corollary 7.1k]{Sch}) says that such a system has a solution iff the following conditions are satisfied for all row vectors $\y,\z \ge \0$: \par
(i) if $\y A + \z B = \0$ then $\y\b + \z\c \ge \0$,\par
(ii) if $\y A + \z B = \0$ and $\y \neq 0$ then $\y\b + \z\c > \0$.\par
For $\b = \c = \0$, condition (i) is always satisfied and condition (ii) is satisfied iff there are no vectors $\y,\z \ge 0$ with $\y A + \z B = \0$ and $\y \neq \0$.
\end{proof}

Using Lemma \ref{Motzkin} we obtain the following criterion of the Cohen-Macaulayness of all symbolic powers of Stanley-Reisner ideals, which is the first step in the proof that simplicial complexes with this property are exactly matroid complexes.
\smallskip

For a subset $F$ of $[n]$ we denote by $\a_F$ the incidence vector of $F$ (which has the $i$-th component equal to 1 if $i \in F$ and 0 else).

\begin{Theorem}\label{all}
Let $\D$ be a pure simplicial complex. Then $I_\D^{(m)}$ is Cohen-Macaulay for all $m \ge 1$ if and only if for every non-Cohen-Macaulay subcomplex $\Gamma$ of $\D$ with $\F(\Gamma) \subseteq \F(\D)$, there exist facets $F_1,...,F_s \in \F(\D)\setminus \F(\Gamma)$ and $G_1,...,G_s \in \F(\Gamma)$ not necessarily different such that 
$$\a_{F_1} + \cdots + \a_{F_s} = \a_{G_1} + \cdots + \a_{G_s}.$$
\end{Theorem}

\begin{proof}
By Theorem \ref{criteria}(iii) we have to check when $L_\G(I^{(m)}) = \emptyset$ for all $m \ge 1$. By definition, 
$$L_\G(I^{(m)}) = \big\{\a \in \NN^n|\ x^\a \in P_F^m\ \text{for}\ F \in \F(\D) \setminus \F(\G)\ \text{and}\ x^\a \not\in P_G^m\ \text{for}\ G \in \F(\G)\big\}.$$ 
Thus, $L_\G(I^{(m)}) = \emptyset$ for all $m \ge 1$ means that the system
\begin{align*}
\sum_{i \not\in F}a_i \ge m &\ \ \big(F \in \F(\D) \setminus \F(\G)\big),\\
\sum_{i \not\in G}a_i < m &\ \ \big(G \in \F(\G)\big),
\end{align*}
has no solution $\a \in \NN^n$ for all $m \ge 1$. 
This condition is equivalent to the condition that the system 
$$\sum_{i \not\in F}a_i > \sum_{i \not\in G}a_i\ \ \ \big(F \in \F(\D) \setminus \F(\G), G \in \F(\G)\big)$$
has no solution $\a \in \NN^n$.
In fact, any solution $\a\in \NN^n$ of the second system will be a solution of the first system for
$m = \min\big\{\sum_{i \not\in F}a_i|\ F \in \F(\D) \setminus \F(\G)\big\}.$
So we have to study when the homogeneous system
\begin{align*}
\sum_{i \not\in G}a_i - \sum_{i \not\in F}a_i < 0 &\ \ \big(F \in \F(\D) \setminus \F(\G), G \in \F(\G)\big),\\
a_i \ge 0 &\ \ ( i = 1,...,n)
\end{align*}
has no solution $\a \in \RR^n$ because any solution in $\RR^n$ can be replaced by a solution in $\QQ^n$, 
which then leads to a solution in $\NN^n$. \par

Let $A$ and $B$ denote the matrices of the coefficients of the inequalities of the first and second line, respectively.
By Lemma \ref{Motzkin}, the above homogeneous system has no solution iff there exist row vectors $\y,\z \ge \0$ such that $\y A + \z B = \0$ and $\y \neq \0$. Let $c_1,...,c_r$ be the non-zero components of the vector $\y$. Since the rows of $A$ are of the form $\a_{\overline G} - \a_{\overline F}$, where $\overline G$ and $\overline F$ denote the complements of $G$ and $F$, and since  $\a_{\overline G} - \a_{\overline F} = \a_F - \a_G$, we have
$$\y A = c_1(\a_{F_1} - \a_{G_1}) + \cdots + c_r(\a_{F_r} - \a_{G_r})$$
for not necessarily different $F_1,...,F_r \in \F(\D) \setminus \F(\G)$ and $G_1,...,G_r \in \F(\G)$.  Since $B$ is the unit matrix, the relation $\y A + \z B = \0$ just means that the monomial $x^{c_1\a_{F_1}}\cdots x^{c_r\a_{F_r}}$  is divided by the monomial $x^{c_1\a_{G_1}}\cdots x^{c_r\a_{G_r}}$.
On the other hand, since $F_1,...,F_r$ and $G_1,...,G_r$ have the same number of elements, $\deg x^{\a_{F_i}} = \deg x^{\a_{G_j}}$ for all $i, j = 1,...,r$. Therefore,
$$\deg x^{c_1\a_{F_1}}\cdots x^{c_r\a_{F_r}} = \deg x^{c_1\a_{G_1}}\cdots x^{c_r\a_{G_r}}.$$
So we must have $x^{c_1\a_{F_1}}\cdots x^{c_r\a_{F_r}} = x^{c_1\a_{G_1}}\cdots x^{c_r\a_{G_r}}$ or, equivalently,
$$c_1\a_{F_1} + \cdots + c_r\a_{F_r} = c_1\a_{G_1} + \cdots + c_r\a_{G_r}.$$
Replacing $c_i\a_{F_i}$ by $\a_{F_i} + \cdots + \a_{F_i}$ and $c_i\a_{G_i}$ by $\a_{G_i} + \cdots + \a_{G_i}$ ($c_i$ times) we may rewrite the above condition as
$$\a_{F_1} + \cdots + \a_{F_s} = \a_{G_1} + \cdots + \a_{G_s}$$
for not necessarily different $F_1,...,F_s \in \F(\D)\setminus \F(\Gamma)$ and $G_1,...,G_s \in \F(\Gamma)$.
Thus, $L_\G(I_\D^{(m)}) = \emptyset$ for all $m \ge 1$ iff this condition is satisfied.
\end{proof}

The condition of Theorem \ref{all} implies that $\D$ is very compact in the following sense. 
Following the terminology of graph theory we call
 an alternating sequence of distinct vertices and facets $v_1,F_1,v_2,F_2,\ldots, v_t,F_t$ a $t$-{\it cycle} of $\D$
if $v_i,v_{i+1}\in F_i$ for all $i = 1,...,t$, where $v_{t+1} = v_1$.

\begin{Corollary} \label{rectangle}
Assume that $I_\D^{(m)}$ is  Cohen-Macaulay for all $m \ge 1$. 
Then every pair $G_1,G_2$ of facets of $\D$ with $|G_1 \cap G_2| \le \dim \D-1$ is contained in a 4-cycle of $\D$ with vertices outside of $G_1 \cap G_2$ and facets containing $G_1 \cap G_2$. Moreover, one of the vertices of the cycle can be chosen arbitrarily in $G_1 \setminus G_2$ or $G_2 \setminus G_1$.
\end{Corollary}

\begin{proof}
By Theorem \ref{all}, there exist facets  $F_1,...,F_s \neq G_1,G_2$ such that 
$$\a_{F_1} + \cdots + \a_{F_s} = c_1\a_{G_1} + c_2\a_{G_2}$$
for some positive integers $c_1, c_2$, $s = c_1+c_2$. By this relation, $G_1 \cup G_2 = F_1 \cup \cdots \cup F_s$ and every facet $F_i$ contains $G_1 \cap G_2$ and vertices of both $G_1 \setminus G_2$ and $G_2 \setminus G_1$.  
\par

Since $|G_1 \cap G_2| \le |G_1|-2$, we can always find two different vertices in $G_1 \setminus G_2$.
Let $u$ be an arbitrary vertex of $G_1 \setminus G_2$. If for every other vertex $v \in G_1 \setminus G_2$, we have $F_i \cap G_2 = F_j \cap G_2$ for all facets $F_i$ containing $u$ and $F_j$ containing $v$, then $F_i \cap G_2 =  F_j \cap G_2$ for all $i,j = 1,...,s$. Since $G_2 \subset F_1 \cup \cdots \cup F_s$, this implies $G_2 \subseteq F_i$ for all $i = 1,...,s$, a contradiction. Therefore, there is another vertex $v$ in $G_1 \setminus G_2$ such that $F_i \cap G_2 \neq F_j \cap G_2$ for some facets $F_i$ containing $u$ and $F_j$ containing $v$.\par

Since $F_i,F_j$ contain $G_1 \cap G_2$,  $F_i \cap (G_2 \setminus G_1)  \neq F_j \cap (G_2 \setminus G_1)$. 
So we can find two different vertices $u' \in F_i \cap (G_2 \setminus G_1)$ and $v' \in F_j \cap (G_2 \setminus G_1)$. 
Clearly, $u,G_1,v,F_j,v',G_2,u',F_i$ form a 4-cycle of $\D$ with vertices outside of $G_1 \cap G_2$ and facets containing $G_1 \cap G_2$.
\end{proof}

By Corollary \ref{rectangle}, every pair of disjoint facets of $\D$ is contained in a 4-cycle of $\D$.
If $\dim \D = 1$, this means that every pair of disjoint edges is contained in a rectangle.
It turns out that this is also a sufficient condition for the Cohen-Macaulayness of all symbolic powers of $I_\D$.

\begin{Corollary}\label{all graph} {\rm \cite[Theorem 2.4]{MT}}
Let $\D$ be a graph. Then $I_\D^{(m)}$ is  Cohen-Macaulay for all $m \ge 1$ if and only if every pair of disjoint edges of $\D$ is contained in a rectangle.
\end{Corollary}

\begin{proof}
We only need to prove the sufficient part. Let $G_1,G_2$ be two disjoint edges of $\D$. Let $F_1,F_2$ be the other edges of a rectangle of $\D$ containing $G_1,G_2$. Obviously, 
$$\a_{F_1} + \a_{F_2} = \a_{G_1} + \a_{G_2}.$$
Hence the condition of Theorem \ref{all} is satisfied.
\end{proof}

It is easy to see that a graph defines a matroid complex if and only if every pair of disjoint edges is contained in a rectangle.
This fact together with Theorem \ref{tight} suggest that there may be a strong relationship between matroid complexes and the Cohen-Macaulayness of all symbolic powers.  In fact, we can prove the following result. This result is also proved independently by Varbaro \cite[Theorem 2.1]{Va}.

\begin{Theorem}\label{all matroids}
$I_\D^{(m)}$ is Cohen-Macaulay for all $m \ge 1$ if and only if $\D$ is a matroid complex. 
\end{Theorem}

\begin{proof}
Assume that $I_\D^{(m)}$ is Cohen-Macaulay for all $m \ge 1$. We will show that if $I$ and $J$ are two faces of $\D$ with $|I \setminus J| = 1$ and $|J\setminus I| = 2$, then there is a vertex $x \in J \setminus I$ such that $I \cup\{x\}$ is a face of $\D$. By \cite[Theorem 39.1]{Sch2}, this implies that $\D$ is a matroid complex.\par

Choose two facets $G_1 \supset I$ and $G_2 \supseteq J$ such that $|G_1 \cap G_2|$ is as large as possible.
If $G_1$ contains a vertex $x \in J \setminus I$, then $I \cup\{x\}$ is a face of $\D$ because it is contained in $G_1$.
Therefore, we may assume that $G_1$ doesn't contain any vertex of $J \setminus I$. Then
$|G_1\cap G_2| \le  |G_2| - |J \setminus I| = \dim \D - 1.$ 
Let $I \setminus J = \{u\}$. If $u \in G_2$, then $I \subset G_2$ and $I \cup\{x\}$ is a face of $\D$ for any $x \in G_2 \setminus I$. 
If $u \not\in G_2$, using Corollary \ref{rectangle} we can find a facet $F \supseteq G_1 \cap G_2$ such that $F$ contains $u$ and a vertex $u' \in G_2 \setminus G_1$. Therefore, $F \supset I$ and $|F \cap G_2| \ge |(G_1 \cap G_2)\cup\{u'\}|= |G_1 \cap G_2| +1$, a contradiction to the choice of $G_1$ and $G_2$. So we have proved the necessary part of the assertion.\par

Conversely, assume that $\D$ is a matroid complex. We will use induction to show that $\D$ satisfies the condition of Theorem \ref{all}. 
If $n = 2$, the assertion is trivial. So we may assume that $n \ge 3$. 
Let $\G$ be an arbitrary non-Cohen-Macaulay subcomplex of $\D$ with $\F(\G) \subseteq \F(\D)$. \par

Let $\D_1$ and $\G_1$ be the subcomplexes of $\D$ and $\G$ generating by the facets not containing $n$, respectively.
Then $\D_1$ is a matroid complex on $[n-1]$ and $\G_1$ is a subcomplex of $\D_1$ with $\F(\G_1) \subseteq \F(\D_1)$. 
By induction we may assume that $\D_1$  satisfies the condition of Theorem \ref{all}. 
If $\G_1$ is not Cohen-Macaulay, there exist facets $F_1,...,F_s \in \F(\D_1)\setminus \F(\Gamma_1)$ and $G_1,...,G_s \in \F(\Gamma_1)$ such that 
$$\a_{F_1} + \cdots + \a_{F_s} = \a_{G_1} + \cdots + \a_{G_s}.$$
Clearly, $F_1,...,F_s \in \F(\D)\setminus \F(\Gamma)$ and $G_1,...,G_s \in \F(\Gamma)$. 
So we may assume that $\G_1$ is Cohen-Macaulay. \par

Let $\D_2 = \lk_\D \{n\}$ and $\G_2 = \lk_\G\{n\}$. 
Then $\D_2$ is also a matroid complex on $[n-1]$ and $\G_2$ is a subcomplex of $\D_2$ with $\F(\G_2) \subseteq \F(\D_2)$. 
By induction we may assume that  $\D_2$  satisfies the condition of Theorem \ref{all}. 
If $\G_2$ is not Cohen-Macaulay, there exist facets $F_1',...,F_s' \in \F(\D_2)\setminus \F(\G_2)$ and $G_1',...,G_s' \in \F(\Gamma_2)$ such that 
$$\a_{F_1'} + \cdots + \a_{F_s'} = \a_{G_1'} + \cdots + \a_{G_s'}.$$
Set $F_i = F_i'\cup\{n\}$, and $G_i = G_i'\cup\{n\}$ for all $i = 1,...,s$. Then $F_1,...,F_s \in \F(\D)\setminus \F(\Gamma)$ and $G_1,...,G_s \in \F(\Gamma)$. Clearly,
$$\a_{F_1} + \cdots + \a_{F_s} = \a_{G_1} + \cdots + \a_{G_s}.$$
So we may assume that $\G_2$ is Cohen-Macaulay. \par

Let $\G_2^*$ be the subcomplex of $\G$ generating by the facets containing $n$.
Then  $\G = \G_1 \cup \G_2^*$. 
Since $\G_2^*$ is a cone over $\G_2$, $\G_2^*$ is Cohen-Macaulay. 
Note that $\G_1 \cap \G_2^* \subseteq \G_2$.
If $\G_1 \cap \G_2^* = \G_2$, there is an exact sequence
$$0 \to k[\G] \to k[\G_1]\oplus k[\G_2^*] \to k[\G_2] \to 0.$$
Since $k[\G_1], k[\G_2^*]$ and $k[\G_2]$ are Cohen-Macaulay
with $\dim k[\G_1] = \dim k[\G_2^*] = \dim k[\G_2] +1$, 
we can conclude that $k[\G]$ is Cohen-Macaulay, which  contradicts the assumption that $\G$ is not Cohen-Macaulay.
So $\G_1 \cap \G_2^*$ is properly contained in  $\G_2$. \par
 
Choose  $G_1 \in \F(\G_1)$ and $G_2  \in \F(\G_2^*)$ such that $G_2 \setminus \{n\} \in \G \setminus (\G_1 \cap \G_2^*)$.   
By the definition of matroids there is a vertex $x \in G_1 \setminus G_2$ such that  
$F = (G_2  \setminus \{n\}) \cup \{x\}$ is a facet of $\D$.  Since $G_2  \setminus \{n\} \not\in  \G_1$, $F \not\in \F(\G)$.
By the proof of  Theorem 3.2,  if the condition of Theorem 3.2 is not satisfied for $\G$, the linear inequality
$\sum_{i \not\in F}a_i > \sum_{i \not\in G_2}a_i$ 
has a solution $\a \in \NN^n$.  From this it follows that $a_n > a_x$.
Since $n$ can be chosen to be any vertex, this implies that the coordinates of $\a$ have no minimum, a contradiction.
So we have proved that $\D$ satisfies the condition of Theorem 3.2.
\end{proof}

Theorem \ref{all matroids} has some interesting consequences. First of all, it implies that the Cohen-Macaulayness of all symbolic powers of Stanley-Reisner ideals doesn't depend on the characteristic of the base field.
\smallskip

Given an integer $d \ge 0$, the $d$-skeleton of a simplicial complex is the set of all faces of dimension  $\le d$.
Obviously, every skeleton of a matroid complex is again a matroid complex.

\begin{Corollary}\label{skeleton}
Let $\D$ be a skeleton of a simplex. Then $I_\D^{(m)}$ is Cohen-Macaulay for all $m \ge 1$.
\end{Corollary}

It is known that for a radical ideal $I \subset S$, $I^m$ is Cohen-Macaulay for all $m \ge 1$ if and only if $I$ is a complete intersection \cite{AV}, \cite{Wa}. This phenomenon doesn't hold for the symbolic powers. For instance, if $\D$ is the $d$-skeleton of a simplex, then $I_\D$ is generated by all squarefree monomials of degree $d+2$, which is not a complete intersection if $d \le n-3$.
\smallskip

Following \cite{Sta} we call a simplicial complex $\D$ a {\it flag complex} if all minimal non-faces consist of two elements.
This is equivalent to say that $I_\D$ is the edge ideal of a simple graph. The Cohen-Macaulayness of symbolic powers of such ideals has been studied recently by
Rinaldo, Terai and Yoshida \cite{RTY}. Using the above results we can easily recover one of their main results.

\begin{Corollary}{\rm \cite[Theorem 3.6]{RTY}}
Let $\D$ be a flag simplex and $\G$ the graph of the minimal nonfaces. 
Then $I_\D^{(m)}$ is Cohen-Macaulay for all $m \ge 1$ if and only if $\G$ is a union of disjoint complete graphs.
\end{Corollary}

\begin{proof}
We note first that $\D$ the clique complex of the graph $\bar \G$ of the nonedges of $\G$. 
By \cite[Theorem 3.3]{KOU}, the clique complex of a graph is a matroid complex if and only if there is a partition of the vertices into stable sets such that every nonedge of the graph is contained in a stable set. A stable set of $\bar \G$ is just a complete graph in $\G$. Therefore, there is a partition of $\G$ into complete graphs such that every edge of $\G$ is contained in such a complete graph.
\end{proof}

\section{Preservation of Cohen-Macaulayness}

Let $\D$ be a simplicial complex. We know by \cite[Corollary 2.5]{MT} and \cite[Theorem 3.6]{RTY} that  if $\dim \D = 1$ or $\D$ is a flag complex and if $I_\D^{(t)}$ is Cohen-Macaulay for some $t \ge 3$, then $I_\D^{(m)}$ is Cohen-Macaulay for all $m \ge 1$. So it is quite natural to ask the following questions:

\begin{Question}
{\rm Is $I_\D^{(m)}$ Cohen-Macaulay if $I_\D^{(m+1)}$ is Cohen-Macaulay?}
\end{Question}

\begin{Question}
{\rm Does there exists a number $t$ depending on $\dim \D$ such that if $I_\D^{(t)}$ is Cohen-Macaulay, 
then $I_\D^{(m)}$ is Cohen-Macaulay for all $m \ge 1$?}
\end{Question}

We don't know any counter-example to both questions.
In the following we will prove some related results on the preservation of Cohen-Macaulayness between different symbolic powers.

\begin{Theorem}\label{down}
$I_\D^{(m)}$ is Cohen-Macaulay if $I_\D^{(t)}$ is Cohen-Macaulay for some $t \ge (m-1)^2+1$.
\end{Theorem}

\begin{proof}
Write $t = r(m-1) + s$ with $1 \le s \le m-1$. Then $r \ge m-1 \ge s$.
Assume for the contrary that $I_\D^{(m)}$ is not Cohen-Macaulay. By Theorem \ref{criteria}(iii), there is a non-Cohen-Macaulay subcomplex $\G$ of $\D$ such that $
L_\G(I_\D^{(m)}) \neq \emptyset$. This means that there is  $\a \in \NN^n$ such that
\begin{align*}
\sum_{i \not\in F}a_i \ge m &\ \ \big(F \in \F(\D) \setminus \F(\G)\big),\\
\sum_{i \not\in G}a_i < m &\ \ \big(G \in \F(\G)\big).
\end{align*}
From this it follows that
\begin{align*}
\sum_{i \not\in F}ra_i & \ge rm \ge r(m-1) + s = t\ \ \big(F \in \F(\D) \setminus \F(\G)\big),\\
\sum_{i \not\in G}ra_i & \le r(m-1) < r(m-1) + s = t\ \ \big(G \in \F(\G)\big).
\end{align*}
Thus, $r\a \in L_\G(I^{(t)})$ so that $L_\G(I^{(t)}) \neq \emptyset$. Therefore, $I_\D^{(t)}$ is not Cohen-Macaulay by Theorem \ref{criteria}(iii).
\end{proof}

For $m =2$ we have $(m-1)^2+1 = 2$. Hence Theorem \ref{down} has the following interesting consequence on the Cohen-Macaulayness of the second symbolic power.

\begin{Corollary}\label{down to 2}
If $I_\D^{(t)}$ is Cohen-Macaulay for some $t \ge 3$, then 
$I_\D^{(2)}$ is Cohen-Macaulay.
\end{Corollary}

Using Theorem \ref{second} we obtain strong conditions on simplicial complexes $\D$ for which $I_\D^{(t)}$ is Cohen-Macaulay for some $t \ge 3$. For instance,  the graph of the one-dimensional faces of $\D$ must have diameter $\le 2$ by Corollary \ref{diameter}.
\smallskip 

For $m = 3$ we have $(m-1)^2+1 = 5$. Therefore, $I_\D^{(3)}$ is Cohen-Macaulay if $I_\D^{(t)}$ is Cohen-Macaulay for some $t \ge 5$.
We don't know any example for which $I_\D^{(4)}$ is Cohen-Macaulay but $I_\D^{(3)}$ is not Cohen-Macaulay. 
\smallskip

The next result shows that there exists a number $t$ depending on $n$ such that if $I_\D^{(t)}$ is Cohen-Macaulay, 
then $I_\D^{(m)}$ is Cohen-Macaulay for all $m \ge 1$.

\begin{Theorem}\label{up}
Let $d = \dim \D$. If $I_\D^{(t)}$ is Cohen-Macaulay for some $t \ge (n-d)^{n+1}$, then 
$I_\D^{(m)}$ is Cohen-Macaulay for all $m \ge 1$.
\end{Theorem}

\begin{proof}
Assume for the contrary that $I_\D^{(m)}$ is non-Cohen-Macaulay for some $m \ge 1$.
By Theorem \ref{criteria}(iii),  there is a non-Cohen-Macaulay subcomplex $\G$ of $\D$ with $\F(\G) \subseteq \F(\D)$ such that
$L_\G(I_\D^{(m)}) \neq \emptyset$. This means that the system
\begin{align*}
\sum_{i \not\in F}a_i & \ge m \ \ \big(F \in \F(\D) \setminus \F(\G)\big),\\
\sum_{i \not\in G}a_i & \le m-1\ \ \big(G \in \F(\G)\big),
\end{align*}
has a solution $\a \in \NN^n$. We now consider the system 
 \begin{align*}
\sum_{i \not\in F}a_i -m & \ge 0\ \ \big(F \in \F(\D) \setminus \F(\G)\big),\\
\sum_{i \not\in G}a_i -m& \le 1\ \ \big(G \in \F(\G)\big),\\
m \ge 0,\; a_i & \ge 0 \ \ (i = 1,...,n).
\end{align*}

The solutions of this system in $\RR^{n+1}$ span a rational polyhedron. Let $\x \in \RR^{n+1}$ be a vertex of this polyhedron. Then $\x$ is the solution of a system $A\x = \b$, where $A$ is an $(n+1)\times (n+1)$ matrix and $\b$ a vector with entries $0,\pm 1$. For the $i$-th components $x_i$ of $\x$ we have $x_i = |\det(A_i)|/|\det(A)|$, where $A_i$ is the matrix obtained from the matrix $(A,\b)$ by deleting the column $i$. Putting $a_i = |\det(A_i)|$ for $i = 1,...,n$, we obtain a solution of the first system of inequalities with $m = |\det(A_{n+1})|$. Since the rows of $A_{n+1}$ have at most $n-d$ non-zero entries which are $\pm 1$, their Euclidean norms are $\le \sqrt{n-d}$. Thus, the Hadamard inequality yields
$$|\det(A_{n+1})| \le \sqrt{(n-d)^{n+1}}.$$

So we may assume that $I_\D^{(m)}$ is non-Cohen-Macaulay for some $m \le \sqrt{(n-d)^{n+1}}$. By Theorem \ref{down}, this implies that $I_\D^{(t)}$ is non-Cohen-Macaulay for $t \ge [\sqrt{(n-d)^{n+1}}-1]^2+1$. Since $(n-d)^{n+1} \ge [\sqrt{(n-d)^{n+1}}-1]^2+1$, this gives a contradiction to the assumption that $I_\D^{(t)}$ is Cohen-Macaulay for some $t \ge (n-d)^{n+1}$.
\end{proof}

One may ask what is the smallest number $t_0$ such that if $I_\D^{(t)}$ is Cohen-Macaulay for some $t \ge t_0$, then $I_\D^{(m)}$ is Cohen-Macaulay for all $m \ge 1$. 
\smallskip

By  \cite[Corollary 2.5]{MT} and \cite[Theorem 3.6]{RTY} we have $t_0 = 3$ if $\dim \D = 1$ or if $\D$ is a flag complex. 
For $\dim \D \ge 2$, we only know that $t_0 \ge 3$. In fact, if we consider the simplicial complex $\D^*$ on $[n+1]$ with
$$\F(\D^*) = \{F \cup \{{n+1}\}|\ F \in \F(\D)\big\},$$
then $\dim \D^* = \dim \D + 1$ and $I_{\D^*}^{(m)}$ is the extension of $I_\D^{(m)}$ in $k[x_1,...,x_n,x_{n+1}]$. Therefore,
$I_{\D^*}^{(m)}$ is Cohen-Macaulay if and only if $I_\D^{(m)}$ is Cohen-Macaulay.

\begin{Remark}
{\rm If $\dim \D = 1$,  one can easily find examples such that $I_\D^{(2)}$ is Cohen-Macaulay but $I_\D^{(m)}$ is not Cohen-Macaulay for all $m \ge 3$. By \cite[Corollary 2.5]{MT}, an instance is a cycle of length 5, which is of diameter 2 but has pairs of disjoint edges not contained in any rectangle. Now, starting from an example in the case $\dim \D = 1$, we can construct a simplicial complex $\D$ of any dimension $\ge 2$ such that $I_\D^{(2)}$ is Cohen-Macaulay but $I_\D^{(m)}$ is not Cohen-Macaulay for all $m \ge 3$.}
\end{Remark}

We have found many examples with $\dim \D = 2$ such that $I_\D^{(2)}$ is Cohen-Macaulay.
In all these cases, we either have $I_\D^{(m)}$ Cohen-Macaulay or not Cohen-Macaulay for all $m \ge 3$.
This suggests that we may have $t_0 = 3$ in the case $\dim \D = 2$. In the following we present a simple example with $\dim \D = 2$ such that $I_\D^{(2)}$ is Cohen-Macaulay but  $I_\D^{(m)}$ is not Cohen-Macaulay for all $m \ge 3$, which is not originated from the case $\dim \D =1$.

\begin{Example}{\rm
Let $\D$ be the simplicial complex with the facets
$$\{1,2,3\},\{1,2,4\},\{1,3,4\},\{2,3,4\},\{3,4,5\}.$$ 

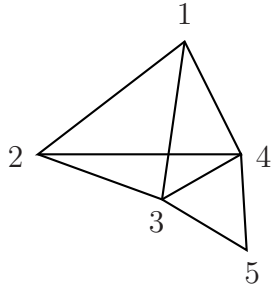
\begin{center}
\psset{unit=1.5cm}
\begin{pspicture}(-2,-1.2)(2,1.4)
\pspolygon(-0.1,1)(-1.4,0)(0.4,0)
\pspolygon(0.4,0)(-0.3,-0.4)(0.45,-0.85)
 \psline(-1.4,0)(-0.3,-0.4)
 \psline(-0.3,-0.4)(-0.1,1)
 \rput(-0.1,1.25){1}
 \rput(-1.6,-0.02){2}
 \rput(-0.35,-0.6){3}
 \rput(0.6,-0.02){4}
 \rput(0.5,-1.05){5}
\end{pspicture}
\end{center}

\noindent We can easily check that $\D$ is a tight complex. By Theorem \ref{tight}, this implies that $I_\D^{(2)}$ is Cohen-Macaulay. To check the Cohen-Macaulayness of $I_\D^{(m)}$, $m \ge 3$, we consider the non-Cohen-Macaulay subcomplex $\G$ with the facets $\{1,2,3\},\{3,4,5\}$. Then $L_\G(I_\D^{(m)}) \neq \emptyset$ because the associated system of linear inequalities:
$$a_1 + a_5 \ge m,\ a_2 + a_5 \ge m,\ a_3+a_5 \ge m,\ a_4+a_5 < m,\ a_1+a_2 <m$$
has at least the solution $\a = (1,1,1,0,m-1)$. By Theorem \ref{criteria}(iii), this implies that $I_\D^{(m)}$ is not Cohen-Macaulay for all $m \ge 3$.
}
\end{Example}

\end{document}